\documentclass[letter,12pt]{article}
\usepackage{amsfonts,amssymb,amsmath,psfrag,graphicx}
\usepackage[bookmarks=true,breaklinks=false]{hyperref}

\newcommand{\comment}[1]{}
\newtheorem{theorem}{Theorem}
\newtheorem{proposition}{Proposition}
\newtheorem{definition}{Definition}

\newlength{\bufQED}

\newcommand{\QED}{\hfill $\!\square$}

\newcommand{\BigOmega}{\mathrm{\Omega}}
\newcommand{\BigTheta}{\mathrm{\Theta}}

\newcommand{\A}{\mathcal{A}}
\newcommand{\e}{\varepsilon}
\renewcommand{\S}{\mathcal{S}}
\newcommand{\C}[1]{\mathcal{C}_{#1}}
\newcommand{\Cbsgs}[1]{\C{#1}^\textrm{\textup{bsgs}}}
\newcommand{\Cbsgsone}[1]{\C{#1}^\textrm{\textup{bsgs1}}}

\newcommand{\Z}{\mathbb{Z}}

\newcommand{\random}{\in_R}
\newcommand{\lead}[1]{\noindent {\bf #1 \ }}

\newcommand{\change}[1]{}

\newcommand{\mystyle}{\setlength{\topsep}{0in}\setlength{\itemsep}{0ex}\setlength{\parskip}{0in}}
\newcounter{step}
{\end{list}}

\title{\vskip2cm Hard Instances of the Constrained Discrete
Logarithm Problem}

\date{}

\author{Ilya Mironov$^1$ \and Anton Mityagin$^{2,4}$ \and
Kobbi Nissim$^{3,4}$\\
}

\begin{document}
\footnotetext[1]{Microsoft Research, Silicon Valley
Campus.}

\footnotetext[2]{Department of Computer Science, University
of California at San Diego, La Jolla, CA 92037.}

\footnotetext[3]{Department of Computer Science, Ben Gurion
University, Beer-Sheva 84105, Israel.}

\footnotetext[4]{The work was done in Microsoft Research,
Silicon Valley Campus.}

\maketitle

\abstract{The discrete logarithm problem (DLP) generalizes
to the constrained DLP, where the secret exponent $x$
belongs to a set known to the attacker. The complexity of
generic algorithms for solving the constrained DLP depends
on the choice of the set. Motivated by cryptographic
applications, we study sets with succinct representation
for which the constrained DLP is hard. We draw on earlier
results due to Erd\"os et~al. and Schnorr, develop
geometric tools such as generalized Menelaus' theorem for
proving lower bounds on the complexity of the constrained
DLP, and construct sets with succinct representation with
provable non-trivial lower bounds.}

\section{Introduction}
\vskip-13cm{\sl An extended abstract of this paper appears
in \emph{7th Algorithmic Number Theory Symposium (ANTS
VII)}, F.~Hess, S.~Pauli, and M.~Pohst, editors, LNCS 4076,
pp. 582--598, 2006. This is the full and corrected
version.}\vskip12cm One of the most important assumptions
in modern cryptography is the hardness of the discrete
logarithm problem (DLP). The scope of this paper is
restricted to groups of prime order $p$, where the DLP is
the problem of computing $x$ given $(g,g^x)$ for $x$ chosen
uniformly at random from $\Z_p$ (see the next section for
notation). In some groups the DLP is believed to have
average complexity of $\BigTheta(\sqrt{p})$ group
operations. The \emph{constrained DLP} is defined as the
problem of computing $x$ given $(g,g^x)$ where $x$ is
chosen uniformly at random from a publicly known set
$S\subseteq\Z_p$.

For the standard DLP there is a well-understood dichotomy
between generic algorithms, which are oblivious to the
underlying group, and group-specific algorithms. By
analogy, we distinguish between generic and group-specific
algorithms for the constrained DLP. In this paper we
concentrate on the former kind, i.e., generic algorithms.
Our main tool for analysis of generic algorithms is the
Shoup-Nechaev generic group
model~\cite{shoup97-generic,Nechaev94}.

The main motivation of our work is the fundamental nature
of the problem and the tantalizing gap that exists between
lower and upper bounds on the constrained DLP. A trivial
generalization of Shoup's proof shows that the DLP
constrained to any set $S\subseteq \Z_p$ has generic
complexity $\BigOmega(\sqrt{|S|})$ group operations. On the
other hand, Schnorr demonstrates that the DLP constrained
to a \emph{random} $S$ of size $\sqrt{p}$ has complexity
$\tilde \BigTheta(\sqrt{p})=\tilde \BigTheta(|S|)$ with
high probability~\cite{Schnorr01-hardcore}. \emph{Explicit}
(de-randomized) constructions or even succinct
representation of small sets with high complexity, or any
complexity better than the square root lower bound were
conspicuously absent.

The importance of improving the square root lower bound for
concrete subsets of $\Z_p$ is implicit
in~\cite{Yacobi98,HS03-sparse,SJ04-hash}, which suggest
exponentiation algorithms that are faster than average for
exponents sampled from certain subsets. These algorithms
either rely on heuristic assumptions of security of the DLP
constrained to their respective sets or use the square root
lower bound to the detriment of their efficiency. For
example, Yacobi proposes to use ``compressible'' exponents
whose binary representation contains repetitive
patterns~\cite{Yacobi98}, which can be exploited by some
algorithms for fast exponentiation. However, without
optimistic assumptions about the complexity of the DLP
constrained to this set the method offers no advantage over
the sliding window exponentiation. Another method of
speeding up exponentiation is to generate an exponent
together with a short addition chain for
it~\cite[Ch.~4.6.3]{knuth-vol2}. Absent reliable methods of
sampling addition chains with uniformly distributed last
elements, this approach depends on the hardness of the DLP
on a non-uniform distribution.

The main technical contribution of our work is the proof
that the DLP constrained to a set $S$, which is chosen from
an easily sampleable family of sets of cardinality
$p^{1/12-\e}$, has complexity $\BigOmega(|S|^{3/5})$ with
probability $1-6p^{-12\e}$. At a higher level of
abstraction we develop combinatorial techniques to bound
the complexity of the constrained DLP, which is a global
property, using the set's local properties. We view our
work as a step towards better understanding the constrained
DLP and possibly designing fast exponentiation algorithms
tuned to work on exponents from ``secure" subsets.


The structure of the paper is as follows. In
Section~\ref{s:generic} we present a number of results
which are known but otherwise scattered in the literature.
In Sections~\ref{s:more} and~\ref{s:beyond} we give new
constructions of sets with provable lower bounds on various
families of algorithms for solving the constrained DLP.

\subsection{Notation}
We use the standard notation for asymptotic growth of
functions, where
\begin{align*}
O(g)&=\{f\colon \mathbb{N}\to\mathbb{R}^+\mid\exists
c,n_0>0\textrm{ s. t. } 0\leq f(n)\leq cg(n)\textrm{ for
all }n>n_0\};\\
\mathrm{\BigOmega}(g)&=\{f\colon
\mathbb{N}\to\mathbb{R}^+\mid\exists c,n_0>0\textrm{ s. t.
} cg(n)\leq f(n)\textrm{ for all
}n>n_0\};\\
\BigTheta(g)&=\{f\colon f=O(g) \textrm{ and }g=O(f)\};\\
\tilde O, \tilde\BigOmega, \tilde \BigTheta&\textrm{ ---
same as }O,\BigOmega,\BigTheta\textrm{ with logarithmic factors ignored;}\\
\Z_p&\textrm{ --- the field of residues modulo prime
$p$;}\\
x\random S&\textrm{ --- $x$ chosen uniformly at random from
$S$.}
\end{align*}

\subsection{Previous work}

Algorithms for solving number-theoretic problems can be
grouped into two main classes: generic attacks, applicable
in any group, and specific attacks designed for particular
groups. The generic attacks on discrete logarithm include
the baby-step giant-step attack~\cite{Shanks71}, Pollard's
rho and lambda algorithms~\cite{Pollard78} as well as their
parallelized versions~\cite{OW99-dlog,Pollard00-kangaroos},
surveyed in~\cite{Teske01-survey}. The specific attacks
have sprouted into a field in their own right, surveyed
in~\cite{SWD96,Odlyzko00}.

A combinatorial view on generic attacks on the DLP was
first introduced by Schnorr~\cite{Schnorr01-hardcore}. He
suggested the concept of the generic DL-complexity of a
subset $S\subseteq \Z_p$ defined as the minimal number of
generic operations required to solve the DLP for any
element of $\{g^x\mid x\in S\}$. He showed that the generic
DL-complexity of random sets of size $m<\sqrt{p}$ is
$m/2+o(1)$. In part our work is an extension of Schnorr's
paper. The combinatorial approach to the DLP was further
advanced by~\cite{CLS03} which gave a characterization of
generic attacks on the entire group of prime order.

Systematically the constrained DLP has been studied for two
special cases: Exponents restricted to an interval and
exponents with low Hamming weight. Pollard's kangaroo
method has complexity proportional to the square root of
the size of the interval~\cite{Pollard00-kangaroos}. The
running time of a simple Las Vegas baby-step giant-step
attack on low-weight exponents is $O(\sqrt t{n/2\choose
t/2})$, where $n$ is the length and $t$ is the weight of
the exponent~\cite{heiman92} (for a deterministic version
see~\cite{Stinson02-bsgs}, which credits it to
Coppersmith). See~\cite{CLP05-bsgs} for cryptanalysis of a
similar scheme in a group of unknown order.


Erd\"os and D.~Newman studied the BSGS-1 complexity (in our
notation) and asked for constructions of sets with a high
(better than a random subset's) BSGS-1
complexity~\cite{EN77-bases}.

\subsection{Generic algorithms}

The generic group model introduced by Shoup and Nechaev
\cite{shoup97-generic,Nechaev94} provides access to a group
$G$ via a random injective mapping $\sigma:G\to\Sigma$,
which \emph{encodes} group elements. The group operation is
implemented as an oracle that on input $\langle \sigma(g),
\sigma(h),\alpha,\beta\rangle $ outputs $\sigma(g^\alpha
h^\beta)$ (for the sake of notation brevity we roll three
group operations, group multiplication, group inversion,
and group exponentiation, in one). Wlog, we restrict the
arguments of the queries issued by algorithms operating in
this model to encodings previously output by the oracle.

The discrete logarithm problem for groups of prime order has a
trivial formalization in the generic group model:
\begin{quote}
Given $p,\sigma(g),\sigma(g^x)$ where $g$ has order $p$ and
$x\random\Z_p$, find $x$.
\end{quote}

The proof sketch of the theorem below, which is essentially
the original one due to Shoup, is reproduced here because
it lays the ground for a systematic study of complexity of
algorithms in the generic group model.

\begin{theorem}[\cite{shoup97-generic}]\label{th:shoup}Let
$\A$ be a probabilistic algorithm and $m$ be the number of
queries made by $\A$. $\A$ solves the discrete logarithm
problem in a group of prime order $p$ with probability
$$
\Pr[\A(p,\sigma(g),\sigma(g^x))=x]<\frac{(m+2)^2}{2p}+\frac1p.
$$
The probability space is $x$, $\A$'s coin tosses, and the
random function $\sigma$.
\end{theorem}
\begin{proof}[sketch] Instead of letting $\A$ interact with a real
oracle, consider the following game played by a simulator.
The simulator keeps track of two lists of equal length
$L_1$ and $L_2$: the list of encodings
$\sigma_1$,\dots,$\sigma_t\in \Sigma$ and the list of
linear polynomials $a_1x+b_1$,\dots,$a_tx+b_t\in \Z_p[x]$.
Initially $L_1$ consists of two elements $\sigma_1$ and
$\sigma_2$, which are the two inputs of $\A$, and $L_2$
consists of $1$ and $x$. When $\A$ issues a query $\langle
\sigma_i,\sigma_j,\alpha,\beta\rangle$, the simulator
fetches the polynomials $a_ix+b_i$ and $a_jx+b_j$ from
$L_2$, computes $a=\alpha a_i+\beta a_j$ and $b=\alpha
b_i+\beta b_j$ and looks up $ax+b$ in $L_2$. If
$ax+b=a_kx+b_k$ for some $k$, the simulator returns
$\sigma_k$ as the answer to the query. Otherwise, the
simulator generates a new element
$\sigma_{t+1}\random\Sigma\setminus L_1$, appends
$\sigma_{t+1}$ to $L_1$ and $ax+b$ to $L_2$, and returns
$\sigma_{t+1}$.

$\A$ terminates by outputting some $y\in \Z_p$. The game
completes as follows:
\begin{enumerate}\mystyle
\item The simulator randomly selects $x^*\random\Z_p$.

\item Compute $a_ix^*+b_i$ for all $i\leq m+2$. If
$a_ix^*+b_i=a_jx^*+b_j$ for some $i\neq j$, the simulator
fails.

\item $\A$ succeeds if and only if $x^*=y$.
\end{enumerate}
Observe that the game played by the simulator is
indistinguishable from the transcript of $\A$'s interaction
with the actual oracle unless the simulator fails in step
2. Since for any two distinct polynomials $a_ix+b_i$ and
$a_jx+b_j$ the probability that $a_ix^*+b_i=a_jx^*+b_j$ is
at most $1/p$, the probability that step 2 fails is at most
$(m+2)^2/{2p}$. Finally, we observe that the probability
that $\A$ wins the game in step 3 is exactly $1/p$, which
completes the proof.\QED
\end{proof}

It follows from the proof that the probability of success
of any probabilistic adaptive algorithm for solving the
discrete logarithm in $\Z_p$ in the generic group model can
be computed given the list of the linear polynomials
induced by its queries. This observation leads us to the
concept of generic complexity defined in the next section.

\section{Generic Complexity}\label{s:generic}

\begin{definition}[Intersection set] For a set of pairs $L\subseteq \Z_p^2$, we define its \emph{intersection
set}
$$I(L)=\{x\in
\Z_p\mid \exists (a,b),(a',b')\in L\textrm{\ s.t.\ }
ax+b=a'x+b'\textrm{\ and\ }(a,b)\neq(a',b')\}.$$
\end{definition}

The set of pairs from the above definition corresponds to
the set of queries asked by the generic algorithm. Its
intersection set is the set of inputs on which the
simulator from the proof of Theorem~\ref{th:shoup} fails.

\begin{definition}[$L$ recognizes an $\alpha$-fraction of $S$] For $L\subseteq \Z_p^2$,
$S\subseteq \Z_p$, and $0<\alpha\leq 1$ we say that $L$
\emph{recognizes an $\alpha$-fraction of} the set $S$ if
$$|S\cap I(L)|\geq\alpha |S|.$$
\end{definition}

\begin{definition}[Generic complexity] The set $S\subseteq \Z_p$
is said to have \emph{generic $\alpha$-complexity} $m$
denoted as $\C{\alpha}(S)$ if $m$ is the smallest
cardinality of a set $L$ recognizing an $\alpha$-fraction
of $S$.
\end{definition}

Our definition of generic complexity is slightly different
from a similar concept of the generic DL-complexity put
forth by Schnorr. We only require that the intersection set
$I(L)$ covers a constant fraction of $S$ rather than the
entire set~\cite{Schnorr01-hardcore}. Our definition better
matches the standard practice of cryptanalysis, when an
attack is considered successful if it succeeds on a
nontrivial fraction of the inputs. Moreover, our bounds
exhibit different scaling behavior as a function of
$\alpha$, and by parametrizing the definition with $\alpha$
we make the dependency explicit.

\begin{proposition}[\cite{Schnorr01-hardcore}]\label{prop:twoineq}For any $S\subseteq \Z_p$ the generic $\alpha$-complexity of
the set $S$ is bounded as
$$
\sqrt{2\alpha|S|}< \C{\alpha}(S)\leq \alpha|S|/2+3.
$$
\end{proposition}
\begin{proof} The lower bound follows from the fact that for any
$L\subseteq \Z_p^2$ the cardinality of the intersection set
is bounded as $|I(L)|<|L|^2/2$. Therefore, in order to
cover at least an $\alpha$-fraction of the set, $|L|^2/2$
must be more than $\alpha|S|$.

The upper bound is attained by the following construction.
If $2m\geq \alpha|S|$ and $\{x_1,\dots,x_{2m}\}\subseteq
S$, then an $\alpha$-fraction of $S$ is recognized by $L$
of size $m+2$ defined as
$$
L=\left\{(0,0),(0,1),(\frac1{x_2-x_1},\frac{x_1}{x_1-x_2}),\dots,(\frac1{x_{2m}-x_{2m-1}},\frac{x_{2m-1}}{x_{2m-1}-x_{2m}})\right\},
$$
since $x_i$ and $x_{i-1}$ are the $x$-coordinates of the
points of intersection of the line
$(\frac1{x_i-x_{i-1}},\frac{x_{i-1}}{x_{i-1}-x_i})$ with
lines $y=0$ and $y=1$ respectively.\QED
\end{proof}

\begin{proposition}\label{prop:C(Z_p)}$\sqrt {2\alpha p}<\C{\alpha}(\Z_p)\leq
2\lceil\sqrt{\alpha p}\rceil$.
\end{proposition}
\begin{proof}The lower bound on $\C{\alpha}(\Z_p)$ is by Proposition~\ref{prop:twoineq}.
The upper bound is given by the set $L=\{(0,i), (1,-\lambda
i)\mid 0\leq i< \lambda\}$, where $\lambda=\lceil
\sqrt{\alpha p}\rceil$. Indeed,
$$
I(L)=\bigcup_{0\leq i,j<\lambda} I(\{(0,i),(1,-\lambda
j)\}=\bigcup_{0\leq i,j<\lambda} \{\lambda j+i\},
$$
which covers $[0,\alpha p)$. \QED
\end{proof}

A tighter (up to a constant factor) bound in the general
case and exact values for $\C{1}(\Z_p)$ for small primes
$p<100$ appear in~\cite{CLS03}.

Since the generic complexity is a monotone property, it
follows that for any set $S\subseteq \Z_p$
$$
\C{\alpha}(S)\leq \min({\alpha}|S|/2+3, 2\lceil\sqrt{\alpha
p}\rceil).
$$

Now we are ready to establish the connection between the generic
complexity of a set and the discrete logarithm problem.

\begin{theorem}Let $S\subseteq \Z_p$, $\A_S$ be a generic algorithm that makes $m<\C{\alpha}(S)$ queries
and outputs a number from $\Z_p$. Then its probability of
success is bounded as\footnote{This statement is stronger
than the one in the proceedings version.}
$$\Pr[\A_S(\sigma(g),\sigma(g^x))=x]<\alpha+\frac1{|S|},$$
where the probability is taken over $\A$'s random tape, the
oracle answers, and $x\random S$. The above bound is tight,
i.e., for any set $S$ there is a generic algorithm whose
query complexity is $\C{\alpha}(S)$ and probability of
success is at least $\alpha+1/|S|$.
\end{theorem}
\begin{proof}[sketch] The proof essentially follows that of
Theorem~\ref{th:shoup}. Let $L$ be a set of pairs
$(a_i,b_i)$ constructed by the simulator and $x^\ast\random
S$ be its choice for $x$. The adversary succeeds in two
cases: either $x^\ast$ belongs to the intersection set of
$L$ or $x^\ast$ is the output of $\A_S$. The first
probability is at most $|I(L)\cap S|/|S|< \alpha$ as long
as $m<\C{\alpha}(S)$, the second probability is exactly
$1/|S|$.

The tightness property follows from the definition of
generic complexity. Let $L$ be the set of pairs of size
$\C{\alpha}(S)$ so that $|S\cap I(L)|\geq \alpha|S|$. Query
the oracle $\langle \sigma(g^x),\sigma(g),a,b\rangle$ for
all pairs $(a,b)\in L$. With probability $|I(L)\cap S|/|S|$
there is a collision that gives away $x$, otherwise make a
guess that succeeds with probability $1/|S\setminus
I(L)|$.\QED
\end{proof}

Notice that the theorem above is unconditional and the
adversary is computationally unbounded. In particular, the
adversary is given full access to $S$ and can design an
$S$-specific algorithm. As long as the algorithm has only
oracle access to the group, $\C{\alpha}(S)$ is a lower
bound on the number $m$ of oracle queries needed by the
algorithm to succeed with probability at least
${\alpha}+1/|S|$.

We know that $\C{\alpha}(S)$ can be negligible compared to
$|S|$ (for instance, according to
Proposition~\ref{prop:C(Z_p)}, when $S=\Z_p$, $|S|=p$ but
its generic complexity is $O(\sqrt p)$). Since the generic
complexity is intimately related to the query complexity of
any discrete logarithm-solving algorithm, we would like to
build sets with higher generic complexity. The next theorem
demonstrates that for a fixed $p$ a random set of size less
than $\sqrt p$ has a near-linear generic complexity.

\begin{theorem}\label{th:random}For a random subset $S\subseteq_R \Z_p$ of size
$p^\e$ for some constant $\e\leq 1/2$ its generic
${\alpha}$-complexity is at least
$$
\C{\alpha}(S)>\frac{\alpha |S|}{\ln p}
$$
with probability $1-1/p$ for large enough $p$.
\end{theorem}
\begin{proof} The proof is by a counting argument.
We shall bound the number of the sets $S$ of size $k=p^\e$
whose $\alpha$-fraction can be recognized by a set $L$ of
size $\delta k$, when $\e\leq 1/2$ and $\delta=\alpha/\ln
p$. Suppose $|L|=\delta k$ and $|I(L)\cap S|\geq \alpha k$,
where $S$ is to be constructed. There are ${p^2\choose
\delta k}$ subsets $L\subseteq \Z^2_p$ of size $\delta k$.
The intersection set $I(L)$ has size at most $(\delta k)^2$
and contains at least $\alpha k$ elements which belong to
$S$. There are thus at most $\binom{p^2}{\delta
k}\binom{(\delta k)^2}{\alpha k}$ distinct possibilities
for these $\alpha k$ elements. The $(1-\alpha)$-fraction of
$S$ can be chosen arbitrarily from $\Z_p$, in ${p\choose
(1-\alpha)k}$ many ways. In total the number of subsets $S$
of generic complexity $\delta k$ and cardinality $k$ is
bounded by ${p^2\choose \delta k}{(\delta k)^2\choose
{\alpha k}}{p\choose (1-\alpha)k}$. Using ${n\choose
k}<\left(ne/k\right)^k$ for any $0<k\leq n$ and
$x^{-x}<1.5$ for any $x>0$ we bound the product as
\begin{multline*}p^{2\delta k+(1-\alpha)k}\delta^{-\delta k+2\alpha
k}k^{-\delta k+2\alpha k-k}\alpha^{-\alpha
k}(1-\alpha)^{-(1-\alpha)k}e^{\delta k+k}<\\
\left[4p^{2\delta+(1-\alpha)}\delta^{2\alpha
}k^{-\delta+2\alpha-1}e^{\delta+1}\right]^k<\left[12p^{2\delta+(1-\alpha)}\delta^{2\alpha
}k^{-\delta+2\alpha-1}\right]^k.
\end{multline*}
We want this number to be less than a $1/p$-fraction of the
number of subsets of $\Z_p$ of size $k$, which is
$1/p{p\choose k}>1/p(p/k)^k$. By taking the $k$th root of
both numbers and substituting $k=p^\e$, we arrive at the
following inequality:
$$
12\delta^{2\alpha}p^{2\delta+(1-\alpha)+\e(-\delta+2\alpha-1)+p^{-\e}}<p^{1-\e}.
$$
Notice that the inequality holds for $\delta <
\alpha(1-2\e)$ if $\e<1/2$ and for $\delta<\alpha/\ln p$ if
$\e=1/2$. When $\e$ is constant and $p$ is large enough,
$\delta=\alpha/\ln p<\alpha(1-2\e)$.\QED
\end{proof}

The bottom line of the theorem we just proved is that hard
sets (where the discrete logarithm is hard to compute using
a generic algorithm) are easy to come by. In fact, almost
any set has high generic complexity (also previously
observed in~\cite{Schnorr01-hardcore}).

In what follows we sharply lower the amount of randomness
that is required to provide any non-trivial guarantee of
generic complexity.

\section{More complexities and lower bounds}\label{s:more}

Many sets of group elements with special properties may be
attacked using a baby-step giant-step method. In this
method the attacker first computes
$g^{c_1}$,\dots,$g^{c_m}$ (giant steps) and then compares
them against $g^{a_1x+b_1}$,\dots,$g^{a_mx+b_m}$ (baby
steps). Any collision between a baby step and a giant step
gives away $x$. We define the complexity of this method
along the lines of the generic complexity from the previous
section.

\begin{definition}[Intersection set-2] For a set of pairs $L\subseteq \Z_p^2$
and a set of points $C\subseteq\Z_p$, we define their
\emph{intersection set} as
$$
I(L,C)=\{x\in \Z_p\mid \exists (a,b)\in L, c\in C\textrm{\
s.t.\ }a\neq 0\textrm{\ and\ }ax+b=c\}.
$$
\end{definition}

\begin{definition}[Baby-step giant-step complexity.] The set $S\subseteq \Z_p$
is said to have the \emph{baby-step giant-step
$\alpha$-complexity} (BSGS complexity for short) $m$
denoted as $\Cbsgs{\alpha}(S)$ if $m$ is the smallest
integer such that there exist $L\subseteq \Z_p^2$ and
$C\subseteq\Z_p$, with $|L|=|C|=m$ and $|I(L,C)\cap
S|\geq\alpha |S|$.
\end{definition}

An important particular case of the baby-step giant-step method is
when all lines defined by $L$ are parallel (i.e., all $a_i=1$).

\begin{definition}[BSGS-1 complexity] The set $S\subseteq \Z_p$
has \emph{BSGS-1 $\alpha$-complexity} $m$ denoted as
$\Cbsgsone\alpha(S)$ if $m$ is the smallest integer such
that there exist $L\subseteq \{1\}\times \Z_p$ and
$C\subseteq\Z_p$, with $|L|=|C|=m$ and $|I(L,C)\cap
S|\geq\alpha|S|$.
\end{definition}

\begin{figure}[htb]
\begin{center}
\psfrag{BSGS}{\scriptsize$I(\{(a_1,b_1),(a_2,b_2),(a_3,b_3)\},\{c_1,c_2,c_3\})$}
\psfrag{BSGS1}{\scriptsize$I(\{(1,b_1),(1,b_2),(1,b_3)\},\{c_1,c_2,c_3\})$}
\psfrag{generic complexity}{\scriptsize$I(\{(a_1,b_1),(a_2,b_2),(a_3,b_3),(a_4,b_4),(a_5,b_5)\})$}
\psfrag{x}{$x$}
\psfrag{c1}{$c_1$}
\psfrag{c2}{$c_2$}
\psfrag{c3}{$c_3$}
\psfrag{x+b1}{$x+b_1$}
\psfrag{x+b2}{$x+b_2$}
\psfrag{x+b3}{$x+b_3$}
\psfrag{a1x+b1}{$a_1x+b_1$}
\psfrag{a2x+b2}{$a_2x+b_2$}
\psfrag{a3x+b3}{$a_3x+b_3$}
\psfrag{a4x+b4}{$a_4x+b_4$}
\psfrag{a5x+b5}{$a_5x+b_5$}
\includegraphics[scale=0.9]{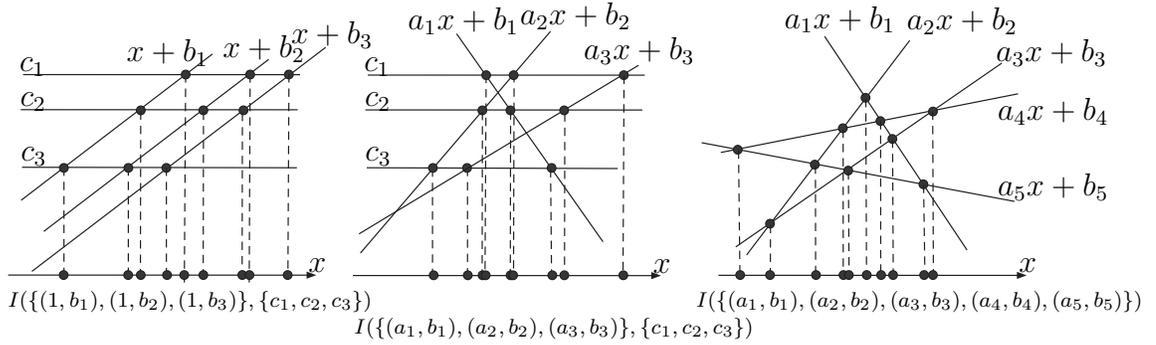}
\caption{\small Intersection sets for BSGS-1, BSGS, and
generic complexities.}\label{fig:intersectionsets}
\end{center}
\end{figure}

Equivalently, $\Cbsgsone\alpha(S)$ is the smallest integer
$n$ such that there exist $X,Y\subseteq \Z_p$, with
$n=|X|=|Y|$ and $|S\cap (X-Y)|\geq\alpha|S|$, where $X-Y$
is the set of pairwise differences between $X$ and $Y$. The
intersection sets from the three definitions of
complexities appear in Fig.~\ref{fig:intersectionsets}.

The problem of computing $\Cbsgsone\alpha(S)$ is
superficially similar to a number of problems in additive
number theory concerned with studying properties of $X-Y$.
However, our goal is fundamentally different since we
require that $X-Y$ cover a large fraction of $S$ rather
than be its exact equal. To the best of our knowledge, the
only paper in the literature directly applicable to
bounding $\Cbsgsone1(S)$ is a 1977 paper by Erd\"os and
Newman~\cite{EN77-bases}. They proved analogues of our
Theorems~\ref{th:random} and~\ref{th:bsgs1} and bounded the
BSGS-1 complexity (in our notation) of the set of small
squares $\{x^2\mid x<\sqrt{p}\}$ to be $\Omega(p^{1/3 -
c/\log\log p})$. They leave open the problem of
constructing sets with a strictly linear BSGS-1 complexity
(without the $1/\log p$ factor).

The BSGS and BSGS-1 complexities provide useful upper
bounds for the generic complexity.

\begin{proposition}\label{prop:bounds}$\frac12\C{\alpha}(S)\leq \Cbsgs{\alpha}(S)\leq \Cbsgsone\alpha(S)$.
\end{proposition}
\begin{proof} Let $C'=\{0\}\times C=\{(0,c)\mid c\in
C\}$. Then $I(L,C)\subseteq I(L\cup C')$, which implies the
first inequality. The second inequality follows from the
fact that any BSGS-1 attack is also a BSGS attack.\QED
\end{proof}

Consider, for example, the baby-step giant-step attacks on
exponents with low Hamming
weight~\cite{heiman92,Stinson02-bsgs}. Define
$S_\lambda=\{x \in \Z_p\mid \nu(x)=\lambda |x|\}$, where
$\nu(x)$ is the number of ones in the binary representation
of $x$. Stinson~\cite{Stinson02-bsgs} estimates the
complexity of the randomized algorithm due to Coppersmith
to yield
$$
\Cbsgsone{1/2}(S_\lambda)=\tilde
O(p^{1/2\log_2(\lambda^{-\lambda}(1-\lambda)^{\lambda-1})}).
$$
For instance, if $\lambda=1/4$, the bound becomes
$\Cbsgsone{1/2}(S_{0.25})=\tilde O(p^{0.406})$.

Following
Propositions~\ref{prop:C(Z_p)},~\ref{prop:bounds}, and
Theorem~\ref{th:random} the BSGS-1 $\alpha$-complexity of a
set of cardinality less than $\sqrt{p}$ lies between
$\sqrt{\alpha|S|/2}$ and $2\sqrt{\alpha p}$, where the
lower bound is trivial and the upper bound is approximated
up to a logarithmic factor by almost any subset of size
$\sqrt{p}$. In this section we construct a set with
succinct representation and a non-trivial BSGS-1
complexity. We start by stating without proof an important
combinatorial lemma known as the Zarankiewicz
problem~\cite{Zarankiewicz51}:
\begin{theorem}\cite[Ch.
IV.2]{Bollobas98-modern} Let $Z(n,s,t)$ be the maximum
number of ones that can be arranged in an $n\times n$
matrix such that there is no all-one $t\times s$ (possibly
disjoint) submatrix. Then
$$
Z(n,s,t)<s^{1/t}n^{2-1/t}.
$$
\end{theorem}
Notice that the asymptotic of the bound on $Z(n,s,t)$
depends on the smallest of the two dimensions of the
prohibited all-one submatrix. It is known that the bound is
tight (up to a constant factor) for $t=2,3$.

Our second combinatorial tool follows from a more general
upper bound due to A.~Naor and Verstra\"ete on the number
of edges in a bipartite graph without cycles of length $2k$
($C_{2k}$-free graph):

\begin{theorem}[\cite{NV05-2kgon}]\label{th:2kgon} The maximum number of edges in a $C_{2k}$-free
$(n,n)$-bipartite graph is less than $2kn^{1+1/k}$.
\end{theorem}

When $k=2$ the two theorems overlap. Indeed, a 0-1 matrix
is also a bipartite graph, where the rows and columns form
the vertex set and the non-zero elements indicate adjacency
of corresponding vertices. In this case an all-one $2\times
2$ submatrix represents a cycle of length 4 in the graph.
Our theorems fully reflect this relationship:
Theorem~\ref{th:bsgs1} can be proved using either the
Zarankiewicz or the Naor-Verstra\"ete bound; its
generalization Theorem~\ref{th:bsgs1general} makes use of
$C_{2k}$-free graphs, while Theorems~\ref{th:bsgs}
and~\ref{th:cbound} apply the Zarankiewicz bound.

\begin{theorem}\label{th:bsgs1}Suppose $S\subseteq \Z_p$ has the property that all
pairwise sums of different elements of $S$ are distinct.
Then
$$
\Cbsgsone\alpha(S)>(\alpha|S|/\sqrt{2})^{2/3}.
$$
\end{theorem}
\begin{proof}Take $X,Y\subseteq\Z_q$, such that $n=|X|=|Y|$ and
$|S\cap (X-Y)|>\alpha|S|$. Consider an $n\times n$ matrix
$M$, whose rows and columns are labeled with elements of
$X$ and $Y$ respectively. For each element $s\in S\cap
(X-Y)$ find one pair $x\in X$ and $y\in Y$ such that
$s=x-y$ and set the $(x,y)$ entry of the matrix to one.
Since $X-Y$ covers at least an $\alpha$-fraction of $S$,
the number of ones in the matrix is at least $\alpha|S|$.

We claim that $M$ does not contain an all-one $2\times 2$
submatrix. Assume the opposite: The submatrix given by
elements $x_1,x_2$ and $y_1,y_2$ has four ones. It follows
that all four $s_{11}=x_1-y_1$, $s_{12}=x_1-y_2$,
$s_{21}=x_2-y_1$, $s_{22}=x_2-y_2\in S$. Then
$s_{11}+s_{22}=s_{12}+s_{21}$, which contradicts the
assumption that all pairwise sums of elements of $S$ are
distinct. Applying the Zarankiewicz bound for the case
$s=t=2$, we prove that
$$
\alpha|S|<Z(n,2,2)<\sqrt{2}n^{2-1/2}=\sqrt{2}n^{3/2},
$$
which implies that
$n=\Cbsgsone\alpha(S)>(\alpha|S|/\sqrt{2})^{2/3}$.\QED
\end{proof}

The sets where sums of pairs of different elements are
distinct are known in combinatorics as weak Sidon sets.
They are closely related to (strong) Sidon sets, also
called $B_2$ sequences, where all pairwise sums (of not
necessarily different elements) are distinct (for a
comprehensive survey see~\cite{Bryant04-survey} that
includes more than 120 bibliographic entries). Explicit
constructions of Sidon subsets of $\{1,\dots,n\}$ due to
Singer and Ruzsa have cardinality at least
$\sqrt{n}-n^{.263}$~\cite{Singer38,BC62-additive,Ruzsa93,BHP01}.

We additionally require that the sums be different modulo
$p$. The size of such sets is bounded from above by
$p^{1/2}+1$~\cite[Theorem 3]{HHO04}. The easiest shortcut
to constructing weak modular Sidon sets is to take a strong
Sidon subset of $\{0\dots \lfloor p/2\rfloor\}$ (see
also~\cite[Ch.~3]{Bryant02-thesis} and~\cite{GS80-golomb}).
Denser Sidon sets may be constructed for primes of the form
$p=q^2+q+1$, where $q$ is also prime~\cite{Singer38}.
Existence of infinitely many such primes is implied by
Schinzel's Hypothesis H and their density follows from the
even stronger Bateman-Horn
conjecture~\cite[A]{Guy04-unsolved}. Interestingly, modular
Sidon sets are useful not only in constructing sets with
high complexity, via Theorem~\ref{th:bsgs1}, but also for
solving the discrete logarithm problem in
$\Z_p$~\cite{CLS03}.

\comment{For a ``book proof'' of the following theorem we
refer the reader to the survey~\cite{Bryant04-survey} that
combines earlier results due to Singer, Erd\"os, Tur\'an,
and others:

\begin{theorem}[\cite{Singer38,ET41-Sidon,Lindstrom69-Sidon,Ruzsa93}]\label{th:sidon} The largest
Sidon subset of $\{1,\dots,n\}$ has between
$\sqrt{n}-n^{.263}$ and $\sqrt{n}+n^{.25}$ elements.
\end{theorem}

Whether the above lower bound can be improved to
$\sqrt{n}+O(1)$ is Erd\"os's \$500-dollar question (problem
C9 in~\cite{guy}).}

\section{Beyond the Basics}\label{s:beyond}

Theorem~\ref{th:bsgs1} can be generalized to make use of
Sidon sets of higher order. First, we prove that if all
$k$-wise sums of elements of $S$ are distinct (counting
permutations of the same $k$-tuple only once), then there
is a bound on the BSGS-1 complexity. Second, we provide a
result that there exist such sets of size
$\BigTheta(p^{1/k})$.

\begin{theorem}\label{th:bsgs1general}Suppose $S\subseteq \Z_p$ is such that all
$k$-wise sums of different elements of $S$ are distinct (excluding
permutation of the summands). Then
$$
\Cbsgsone\alpha(S)>(\alpha|S|/(2k))^{k/(k+1)}.
$$
\end{theorem}
\begin{proof} Take $X,Y\subseteq\Z_p$, such that $\Cbsgsone\alpha(S)=|X|=|Y|=n$ and
$|S\cap (X-Y)|>\alpha|S|$. Instead of the matrix as in
Theorem~\ref{th:bsgs1}, consider a bipartite graph $G(X,Y)$, where
there is an edge $(x,y)$ if and only if $x-y\in S$ (keep only one
edge per element of $S$).

We claim that there are no $2k$-cycles (without repetitive
edges) in the bipartite graph $G$. Assume the opposite:
There is a cycle $(x_1$,$y_1$,\dots,
$x_k$,$y_k$,$x_1$,$y_1)$. Consider two sums:
$(x_1-y_1)+(x_2-y_2)+\dots+(x_k-y_k)$ and
$(x_2-y_1)+(x_3-y_2)+\dots+(x_k-y_{k-1})+(x_1-y_k)$. Not
only are the two sums equal, they also consist of $k$
elements of $S$ each, and these elements are all distinct
(as every element of $S$ appears as an edge of $G$ at most
once). A contradiction is found.

The number of edges in an  $(n,n)$-bipartite graph without
$2k$-cycles is less than $2kn^{1+1/k}$
(Theorem~\ref{th:2kgon}). Therefore $\alpha|S| <
2kn^{1+1/k}$, and
$n=\Cbsgsone\alpha(S)>(\alpha|S|/(2k))^{k/(k+1)}$. \QED
\end{proof}

Bose and Chowla give a construction for subsets of
$\{1,\dots,q^k\}$ of prime size $q$ whose $k$-wise sums are
distinct (in integers, not modulo
$p$)~\cite{BC62-additive}. By choosing the largest prime
$q$ less than $p^{1/k}$ (which, for large $p$ is more than
$p^{1/k}-p^{0.525/k}$~\cite{BHP01}) an interval of length
$q^k/k$ with a $1/k$ proportion of the set's elements, we
guarantee that all $k$-sums are distinct in $\Z_p$ as well.
Unfortunately, \cite{BC62-additive} does not provide an
efficient sampling algorithm.

Along the lines of Theorem~\ref{th:bsgs1} we prove that
other verifiable criteria imply non-trivial bounds on the
BSGS and generic complexity.

\begin{theorem}\label{th:bsgs}Suppose $S\subseteq \Z_p$
is such that for any distinct $x_1$,$x_2$,$y_1$,
$y_2$,$z_1$,$z_2\in S$:
\begin{equation}\label{eq:det2x2}
\det\left(\begin{matrix}x_1-y_1 &
x_2-y_2\\
y_1-z_1 & y_2-z_2\end{matrix}\right)\neq 0.
\end{equation}
Then
$$
\Cbsgs{\alpha}(S)>(\alpha|S|/\sqrt{3})^{2/3}.
$$
\end{theorem}
\begin{proof} Take $L\subseteq \Z_p^2$ and
$C\subseteq\Z_p$, such that $|L|=|C|=n$ and $|I(L,C)\cap
S|>\alpha|S|$. As in Theorem~\ref{th:bsgs1} consider the
$n\times n$ matrix $M$ whose rows and columns are labeled
with elements of $L$ and $C$ respectively. For each element
$s\in S\cap I(L,C)$ set one entry in row $(a,b)$ and column
$c$ to one, where $s=(c-b)/a$. Thus, the total number of
ones in the matrix is exactly $m=|I(L,C)\cap S|$. If there
is a $2\times 3$ all-one submatrix in $M$, then
property~(\ref{eq:det2x2}) does not hold (three parallel
lines divide two other lines proportionally). The
Zarankiewicz bound implies that
$$
\alpha|S|<Z(n,3,2)<\sqrt{3}n^{2-1/2}=\sqrt{3}n^{3/2}.
$$
Hence $\Cbsgs{\alpha}(S)=n\geq (\alpha|S|/\sqrt{3})^{2/3}$.\QED
\end{proof}

Constructing a large subset of $\Z_p$ with short
description satisfying the condition of the previous
theorem is a difficult problem. Fortunately, the
probability that a random 6-tuple of $\Z_p$ elements fails
to satisfy~(\ref{eq:det2x2}) is $2/p$~\cite{Schwartz80}.
This observation motivates the following definition:
\begin{definition}[$\S(N,k)$ family] Let $\S(N,k)=\{x_1,\dots,x_N\}$ be
a family of subsets of $\Z_p$, where $x_1,\dots,x_N\colon
\mathcal{K}\mapsto\Z_p$ are $k$-wise independent random
variables ($\mathcal{K}$ is the probability space).
\end{definition}

Properties of $\S(N,k)$ are established in the following
proposition:
\begin{proposition}\label{prop:Snk}
\begin{enumerate}\mystyle
\item $\S(N,k)$ can be defined over $\mathcal{K}=\Z_p^k$.

\item For $k>1$,
$\Pr_{S\random \S(N,k)}[|S|\neq N]<N^2/p$.

\item If $h\in \Z[y_1,\dots,y_k]$ and $d=\deg(h)>0$, then
$$
\Pr_{S\in \S(N,k)}[\exists\textrm{ distinct }
z_1,\dots,z_k\in S\textrm{ with }
h(z_1,\dots,z_k)=0]<N^kd/p.
$$
\end{enumerate}
\end{proposition}

\begin{proof}
1. To construct $\S(N,k)$ we use a well-known $k$-universal
family of functions (following~\cite{CW77}). Let the
probability space be $\mathcal{K}=\Z_p^k$ and
$f_a(x)=a_{k-1}x^{k-1}+\dots+a_0$ for
$a=(a_0,\dots,a_{k-1})\in \mathcal{K}$. Define the random
variables $x_i=f_a(i)\colon \mathcal{K}\to\Z_p$ for $1\leq
i\leq N$. We claim that the variables $x_1,\dots,x_N$ are
$k$-wise independent. This follows from the system
$f_a(i_1)=y_1$,\dots, $f_a(i_k)=y_k$ having a unique
solution $a\in\mathcal{K}$ for any distinct
$i_1,\dots,i_k\in \{1,\dots,N\}$ and
$y_1,\dots,y_k\in\Z_p$. Notice that any $S\in\S(N,k)$ can
be easily enumerated and sampled from.

2. Let $I_{ij}$ be the indicator variable, which is equal
to 1 when $x_i=x_j$ and 0 otherwise. The cardinality of
$S=\{x_1,\dots,x_N\}$ is at least $N-\sum_{i<j}I_{ij}$.
Since $x_i$ and $x_j$ are independent for all $i\neq j$,
$E[I_{ij}]=1/p$. By linearity of expectation, the expected
value $E[\sum_{i<j}I_{ij}]<N^2/p$. By Markov's inequality
$\Pr[|S|\neq N]=\Pr[\sum_{i<j}I_{ij}\geq 1]<N^2/p$.

3. Let $I_{i_1,\dots,i_k}$ for all distinct $1\leq
i_1,\dots,i_k\leq N$ be the indicator variable that is 1 if
and only if $h(x_{i_1}$,\dots,$x_{i_k})=0$. By independence
of the variables and~\cite{Schwartz80}
$E[I_{i_1,\dots,i_k}]\leq 2/p$, which by linearity of
expectation and Markov's inequality implies that
$\Pr_S[\exists\textrm{ distinct } x_1,\dots,x_k\in
S\textrm{ with } h(x_1,\dots,x_k)=0]\leq
\Pr_S[\sum_{i_1,\dots,i_k} I_{i_1,\dots,i_k}\geq 1]<
N^kd/p.$ \QED
\end{proof}

It follows that a randomly chosen set from the family
$\S(p^{1/6-\e},6)$ has size $p^{1/6-\e}$ with probability
at least $1-p^{-2/3}$ and satisfies the condition of
Theorem~\ref{th:bsgs} with probability at least
$1-2p^{-6\e}$.

\begin{figure}[htb]
\begin{center}
\psfrag{f1}{$x_1+y_2=x_2+y_1$}
\psfrag{f2}{$\left|\begin{matrix}x_1-y_1 &
x_2-y_2\\
y_1-z_1 & y_2-z_2\end{matrix}\right|=0$}
\psfrag{f3}{\scriptsize{$\left|\begin{matrix}
x_1 - y_1 & x_1 - z_1 & z_1(x_1 - y_1) & y_1(x_1 - z_1) \\
x_2 - y_2 & x_2 - z_2 & z_2(x_2 - y_2) & y_2(x_2 - z_2)  \\
x_3 - y_3 & x_3 - z_3 & z_3(x_3 - y_3) & y_3(x_3 - z_3)  \\
x_4 - y_4 & x_4 - z_4 & z_4(x_4 - y_4) & y_4(x_4 - z_4)
\end{matrix}\right|=0$}} 
\psfrag{f4}{$x_1+y_2+z_3=x_3+y_1+z_2$}
\psfrag{x}{$x$}
\psfrag{x1}{$x_1$}
\psfrag{y1}{$y_1$}
\psfrag{z1}{$z_1$}
\psfrag{x2}{$x_2$}
\psfrag{y2}{$y_2$}
\psfrag{z2}{$z_2$}
\psfrag{x3}{$x_3$}
\psfrag{z3}{$z_3$}
\psfrag{x11}{$x_1$}
\psfrag{x12}{$x_2$}
\psfrag{x13}{$x_3$}
\psfrag{x14}{$x_4$}
\psfrag{x21}{$y_1$}
\psfrag{x22}{$y_2$}
\psfrag{x23}{$y_3$}
\psfrag{x24}{$y_4$}
\psfrag{x31}{$z_1$}
\psfrag{x32}{$z_2$}
\psfrag{x33}{$z_3$}
\psfrag{x34}{$z_4$}
\psfrag{l1}{$l_1$}
\psfrag{l2}{$l_2$}
\psfrag{l3}{$l_3$}
\psfrag{l4}{$l_4$}
\psfrag{ll1}{$l_x$}
\psfrag{ll2}{$l_y$}
\psfrag{ll3}{$l_z$}
\includegraphics[scale=0.9]{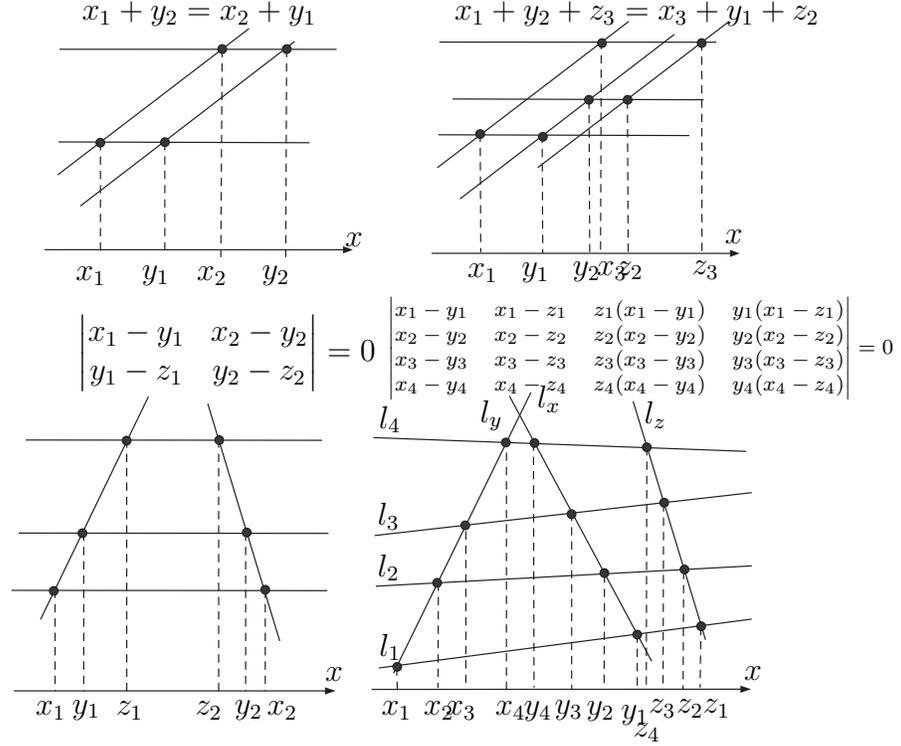}
\caption{\small ``Prohibited'' configurations
(Theorems~\ref{th:bsgs1},~\ref{th:bsgs1general}
and~\ref{th:bsgs},
Proposition~\ref{prop:12}).}\label{fig:configurations}
\end{center}
\end{figure}

To apply a similar argument to the all-powerful generic
complexity, we may show that for small constants $m_1$ and
$m_2$, the projections on the $x$ axis of the intersection
points of an irregular $m_1$ by $m_2$ grid (in which lines
need not be parallel) satisfy a certain relationship. Next,
a set $S$, where any $m_1m_2$-tuple avoids this
relationship, is to be constructed.

Let us see first why this argument works for some values of
$m_1$ and $m_2$, and then improve the parameters. Let
$m_1=4$ and $m_2=5$. There are $9$ lines that can be
described using 18 parameters. On the other hand, there are
20 points that form the intersection set of these lines.
Each of the 20 intersection points imposes a linear
equation on the parameters, and hence the system is
overdetermined (even if we exclude linearly dependent
equations). In particular, this implies that the
probability that a random 20-tuple of elements of $\Z_p$ is
coverable by a $4\times 5$ grid is negligible. We refine
this argument in the following proposition.

\begin{proposition}[Bipartite Menelaus' theorem]\label{prop:12}
Consider seven lines $l_{x,y,z},l_{1,2,3,4}$ forming an irregular
grid, and their twelve intersection points. Let $x_i$,$y_i$,$z_i$ be
projections on the $x$ axis of the intersection points of $l_i$ with
lines $l_x,l_y,l_z$. Then the following holds:
\begin{equation}
\det\left(\begin{matrix}
x_1 - y_1 & x_1 - z_1 & z_1(x_1 - y_1) & y_1(x_1 - z_1) \\
x_2 - y_2 & x_2 - z_2 & z_2(x_2 - y_2) & y_2(x_2 - z_2)  \\
x_3 - y_3 & x_3 - z_3 & z_3(x_3 - y_3) & y_3(x_3 - z_3)  \\
x_4 - y_4 & x_4 - z_4 & z_4(x_4 - y_4) & y_4(x_4 - z_4)
\end{matrix}\right)=0.\label{eq:twelvedet}
\end{equation}
\comment{
\begin{align}
0=\phantom{+}
 &(x_1-z_1) (y_1-y_2) (x_2-z_2) (x_3-x_4) (y_3-z_3) (y_4-z_4) -\notag\\
-&(x_1-z_1) (y_1-y_3) (x_2-x_4) (y_2-z_2) (x_3-z_3) (y_4-z_4) +\notag\\
+&(x_1-z_1) (y_1-y_4) (x_2-x_3) (y_2-z_2) (y_3-z_3) (x_4-z_4) +\notag\\
+&(x_1-x_2) (y_1-z_1) (y_2-z_2) (x_3-z_3) (y_3-y_4) (x_4-z_4) -\notag\\
-&(x_1-x_3) (y_1-z_1) (x_2-z_2) (y_2-y_4) (y_3-z_3) (x_4-z_4) +\notag\\
+&(x_1-x_4) (y_1-z_1) (x_2-z_2) (y_2-y_3) (x_3-z_3) (y_4-z_4).
\end{align}}
\end{proposition}
\begin{proof}Denote the $4\times 4$ matrix
in~(\ref{eq:twelvedet}) by $M$. Observe that if any of the seven
lines is vertical, (\ref{eq:twelvedet}) follows immediately.
Indeed, if $l_y=\{x=\mathrm{const}\}$, then $y_1=y_2=y_3=y_4$ and
the second and the fourth columns of $M$ are linearly dependent.
Moreover, $\det M$ is invariant under permutations of $l_x$,$l_y$,
and $l_z$, which takes care of vertical $l_y$ or $l_z$. If $l_i$
is vertical for some $1\leq i\leq 4$, then the $i$th row of $M$ is
all-zero, and $\det M = 0$.

If none of the lines is vertical, we can write down
equations for all of them in the Cartesian coordinates. Let
$l_{x,y,z}=\{a_{x,y,z}x + b_{x,y,z}\}$ and
$l_i=\{c_ix+d_i\}$ for $1\leq i\leq 4$. Each intersection
point imposes an equation on the parameters of the two
lines incident with it, a total of 12 equations in 14
unknowns.\comment{For instance, the intersection point of
$l_x$ and $l_1$ projects to $x_1$, which implies that
$$
a_xx_1+b_x=c_1x_1+d_1.
$$
There are 12 equations in 14 unknowns if we try to solve
the system for the lines' parameters given
$x_1,x_2,\dots,z_4$.} However, the system always has a
trivial solution, when all lines are equal. Rewrite the
system using new variables: $\tilde a_y=a_y-a_x$, $\tilde
b_y=b_y-b_x$, $\tilde a_z=a_z-a_x$, $\tilde b_z=b_z-b_x$,
$\tilde c_1=c_1-a_x$, $\tilde d_1=d_1-b_x$, etc. The result
is a homogenous system of 12 linear equations in 12 new
variables. It has a non-zero solution if and only if its
matrix is singular (only non-zero elements are shown):
$$
M'=\left(\begin{array}{cccccccccccc}
    \cdot&    \cdot&    \cdot&    \cdot&    \cdot&    \cdot& -x_1&   -1&    \cdot&    \cdot&    \cdot&    \cdot\\
    \cdot&    \cdot&    \cdot&    \cdot&    \cdot&    \cdot&    \cdot&    \cdot& -y_1&   -1&    \cdot&    \cdot\\
    \cdot&    \cdot&    \cdot&    \cdot&    \cdot&    \cdot&    \cdot&    \cdot&    \cdot&    \cdot& -z_1&   -1\\
  x_2&    1&    \cdot&    \cdot&    \cdot&    \cdot& -x_2&   -1&    \cdot&    \cdot&    \cdot&    \cdot\\
  y_2&    1&    \cdot&    \cdot&    \cdot&    \cdot&    \cdot&    \cdot& -y_2&   -1&    \cdot&    \cdot\\
  z_2&    1&    \cdot&    \cdot&    \cdot&    \cdot&    \cdot&    \cdot&    \cdot&    \cdot& -z_2&   -1\\
    \cdot&    \cdot&  x_3&    1&    \cdot&    \cdot& -x_3&   -1&    \cdot&    \cdot&    \cdot&    \cdot\\
    \cdot&    \cdot&  y_3&    1&    \cdot&    \cdot&    \cdot&    \cdot& -y_3&   -1&    \cdot&    \cdot\\
    \cdot&    \cdot&  z_3&    1&    \cdot&    \cdot&    \cdot&    \cdot&    \cdot&    \cdot& -z_3&   -1\\
    \cdot&    \cdot&    \cdot&    \cdot&  x_4&    1& -x_4&   -1&    \cdot&    \cdot&    \cdot&    \cdot\\
    \cdot&    \cdot&    \cdot&    \cdot&  y_4&    1&    \cdot&    \cdot& -y_4&   -1&    \cdot&    \cdot\\
    \cdot&    \cdot&    \cdot&    \cdot&  z_4&    1&    \cdot&    \cdot&    \cdot&    \cdot& -z_4&   -1\\
\end{array}\right).
$$
One can verify that $\det(M)=\det(M')$.\QED
\end{proof}

In the full version of the paper we give a geometric proof
of the proposition, deriving~(\ref{eq:twelvedet}) directly,
and explain the connection with classic Menelaus' theorem.
We also give an alternative statement of the theorem, which
puts it in the realm of projective geometry.

Proposition~\ref{prop:12} is the ``minimal'' condition that
holds for the $x$-coordinates of the intersection points of
two sets of lines in general position. Indeed, it follows
from the proof that for any assignment of distinct values
to the eleven variables $x_{1,2,3,4}$, $y_{1,2,3,4}$,
$z_{1,2,3}$ there is a configuration of lines whose
intersection points project to those variables. Other
configurations with as many or fewer intersection points do
not produce any conditions either. For example, six lines
intersecting two lines can project to any collection of
twelve points.


All geometric arguments
(Theorems~\ref{th:bsgs1},~\ref{th:bsgs1general}, and~\ref{th:bsgs},
Proposition~\ref{prop:12}) are illustrated
in~Fig.~\ref{fig:configurations}.

\begin{theorem}\label{th:cbound}If $S$ is chosen from $\S(p^{1/12-\e},12)$, then
with probability at least $1-6p^{-12\e}$
$$
\C\alpha(S)>(\alpha|S|/\sqrt[3]{4})^{3/5}.
$$
\end{theorem}
\begin{proof}Consider the set of lines $L\subseteq\Z_p^2$ such that
$\C\alpha(S)=|I(L)\cap S|$ and $n=|L|$. As in
Theorem~\ref{th:bsgs}, we apply the Zarankiewicz bound to
the $n\times n$ matrix, only now both the rows and the
columns are labeled with elements of the set $L$.
Similarly, only one occurrence of an element of $S$ as the
$x$-coordinate of the intersection of two distinct lines is
recorded in the matrix.

According to Proposition~\ref{prop:Snk}, $S$ avoids
solutions to the equation~(\ref{eq:twelvedet}), whose
left-hand side is a multivariate polynomial of total degree
6, with probability greater than $1-6p^{-12\e}$. Therefore
the probability that there exist 12 points in $S$ that can
be the intersection set of two groups of lines consisting
of 3 and 4 lines respectively is less than $6p^{-12\e}$.
Finally, as before,
$\alpha|S|<Z(n,4,3)<4^{1/3}n^{2-1/3}=\sqrt[3]{4}n^{5/3}=\sqrt[3]{4}\C\alpha(S)^{5/3}.$\QED
\end{proof}

Unlike the proofs of
Theorems~\ref{th:bsgs1},~\ref{th:bsgs1general},
and~\ref{th:bsgs}, where the classes in which the lines are
grouped arise naturally, the use of bipartite Menelaus'
theorem in the analysis of generic complexity above might
appear less motivated. In fact, classic Menelaus' theorem
imposes a simple condition (a cubic equation) on the
intersection set of four lines. It is the second step of
the argument, where we translate absence of a certain
submatrix (subgraph) into sparseness of the entire matrix,
which becomes problematic: Unless $H$ is bipartite,
$H$-free graphs on $n$ vertices may have as many as
$\BigTheta(n^2)$ edges according to the celebrated Tur\'an
theorem~\cite[Ch. IV.2]{Bollobas98-modern}.

\begin{figure}[htb]
\begin{center}
\psfrag{C(S)}{$\C1,\Cbsgs1,\Cbsgsone1(S)$}
\psfrag{sqrtp}{$\sqrt{p}$} \psfrag{p}{$p$}
\psfrag{p13}{$\sqrt[3]{p}$} \psfrag{p14}{$\sqrt[4]{p}$}
\psfrag{p15}{$\sqrt[5]{p}$} \psfrag{p16}{$\sqrt[6]{p}$}
\psfrag{p19}{$\sqrt[9]{p}$} \psfrag{p112}{$\sqrt[12]{p}$}
\psfrag{p120}{$\sqrt[20]{p}$} \psfrag{S}{$|S|$}
\psfrag{log}{log} \psfrag{bsgs1}{Theorem~\ref{th:bsgs1}
($\Cbsgsone1$)} \psfrag{bsgs}{Theorem~\ref{th:bsgs}
($\Cbsgs1$)} \psfrag{cbound}{Theorem~\ref{th:cbound}
($\C1$)} \psfrag{lemma1}{Prop.~\ref{prop:C(Z_p)}}
\psfrag{lemma2}{Prop.~\ref{prop:twoineq}}
\psfrag{random}{Theorem~\ref{th:random}}
\includegraphics[scale=0.8]{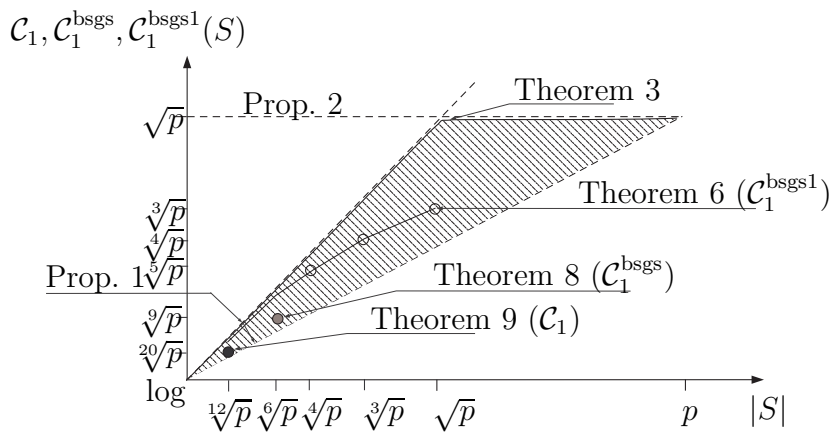}
\caption{\small Generic complexities and bounds (in
logscale). Propositions~\ref{prop:C(Z_p)}
and~\ref{prop:twoineq} bound the triangle that contains all
possible values for generic complexity. The theorems point
to lower bounds provable for complexities of their
respective constructions.\label{fig:allbounds}}
\end{center}
\end{figure}

\section{Conclusion}
In this paper we develop a theory of lower bounds in the
generic group model on the discrete logarithm problem
constrained to a subset $S\subseteq \Z_p$ known to the
attacker (constrained DLP). We give a first construction of
a set with succinct description whose generic complexity is
more than the square root of its size
(Theorem~\ref{th:cbound}). There exists an apparent gap
between our construction ($|S|=p^{1/12}$ and
$\C1(S)=|S|^{3/5}$) and a random set of size $p^{1/2}$
whose complexity is almost linear in its size. Bridging
this gap constitutes an interesting open problem whose
solution would shed some light on the intrinsic difficulty
of the discrete logarithm problem. We also define
restricted versions of the generic complexity that capture
the complexity of baby-step-giant-step algorithms. We give
an explicit, deterministic construction of a collection of
sets, whose complexity in respect to the weakest family of
baby-step-giant-step algorithms becomes near-optimal as
their size decreases (Theorem~\ref{th:bsgs1general}).
Various bounds and constructions are put together in
Fig.~\ref{fig:allbounds}.

\textbf{Acknowledgments.} The authors thank Constantin
Shramov for his crucial contribution to the geometric proof
of Theorem~\ref{th:bipartite}, and Fan Chung, Kevin
O'Bryant, Imre Ruzsa, and the anonymous reviewer of ANTS
VII for their advice and helpful comments.



\appendix

\newcommand{\llog}[1]{\log_2{#1}}

\section{Bipartite Menelaus' Theorem}\label{s:12points}

\begin{figure}[htb]
\begin{center}
\psfrag{x11}{$A_1$}
\psfrag{x12}{$A_2$}
\psfrag{x13}{$A_3$}
\psfrag{x14}{$A_4$}
\psfrag{x21}{$B_1$}
\psfrag{x22}{$B_2$}
\psfrag{x23}{$B_3$}
\psfrag{x24}{$B_4$}
\psfrag{x31}{$C_1$}
\psfrag{x32}{$C_2$}
\psfrag{x33}{$C_3$}
\psfrag{x34}{$C_4$}
\psfrag{l1}{$l_1$}
\psfrag{l2}{$l_2$}
\psfrag{l3}{$l_3$}
\psfrag{l4}{$l_4$}
\psfrag{ll1}{$l_a$}
\psfrag{ll2}{$l_b$}
\psfrag{ll3}{$l_c$}
\psfrag{ll4}{$l_d$}
\psfrag{A}{$A$}
\psfrag{B}{$B$}
\psfrag{C}{$C$}
\psfrag{D}{$D$}
\psfrag{E}{$E$}
\psfrag{F}{$F$}
\psfrag{parta}{a)}
\psfrag{partb}{b)}
\includegraphics{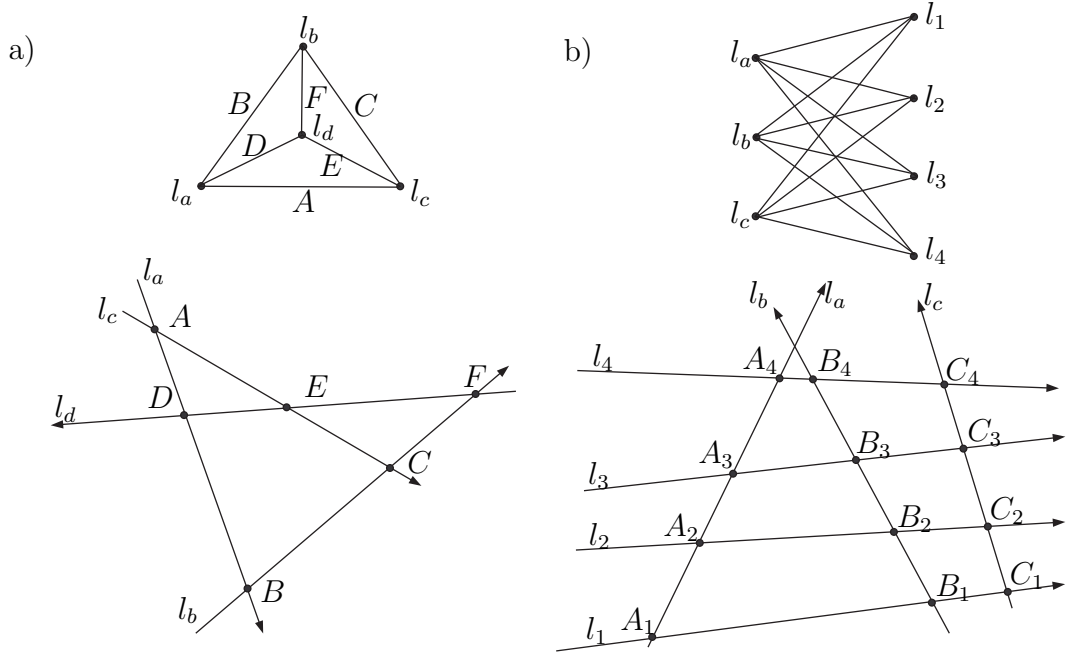}
\caption{Classic and bipartite Menelaus' theorems and their
intersection graphs.}\label{fig:grid}
\end{center}
\end{figure}

Let us first recap the classic Menelaus theorem.
\begin{theorem}[Classic Menelaus]\label{th:classic}Consider four directed lines intersecting at six
points (see Fig.~\ref{fig:grid}a). Then
\begin{equation}\label{eq:menelaus}
AD\cdot BF \cdot CE = BD \cdot CF \cdot AE,
\end{equation}
where the segments' lengths are signed, i.e., positive if their
direction agrees with that of the line they belong to and negative
otherwise.
\end{theorem}

Less known is that (\ref{eq:menelaus}) is equivalent to the
following:
\begin{equation}\label{eq:menelausdet}
\det\left(\begin{matrix}
 DB & AD & AB\\
 BC & FC & BF\\
 ED & DF & FE\\
\end{matrix}\right) = 0.
\end{equation}

One way to interpret the theorem is to add an $x$-axis to the
drawing and consider projections of the intersection points onto
this axis. Since the ratios of signed collinear segments are
invariant under orthogonal projection, (\ref{eq:menelaus}) implies
that
\begin{equation}\label{eq:menelaus-coord}
(x_A-x_D)\cdot (x_B-x_F) \cdot (x_C-x_E) = (x_B-x_D) \cdot (x_C-x_F)
\cdot (x_A-x_E),
\end{equation}
where $x_A$ is the $x$-coordinate of $A$, etc.

We may reverse the Menelaus theorem and ask whether a given
six-tuple can be the projection of the intersection points of four
distinct lines. It is easy to check that (\ref{eq:menelaus-coord})
is not only necessary but is also a sufficient condition for such
four lines to exist. \comment{
\newcommand{\menelaus}[4]{$A_#3A_#4\cdot A_#1L_{#1#4}\cdot B_#1B_#4\cdot B_#2L_{#2#4} \cdot C_#2C_#4\cdot C_#3L_{#3#4}=\
                           A_#1A_#4\cdot B_#1L_{#1#4}\cdot B_#2B_#4\cdot C_#2L_{#2#4} \cdot C_#3C_#4\cdot A_#1L_{#3#4}$}
\menelaus1234
\menelaus2341
\menelaus3412
\menelaus4123}

A natural extension of the classic Menelaus theorem is to consider
other configurations (combinations of lines and their intersection
points). Generalized Menelaus theorem (a line crossing a polygon)
that corresponds to the wheel graph is well known. Below we prove
the smallest possible ``bipartite'' Menelaus theorem.

\begin{theorem}[Bipartite Menelaus]\label{th:bipartite}
Consider seven directed lines $l_{a,b,c},l_{1,2,3,4}$ forming an
irregular grid (see Fig.~\ref{fig:grid}b) and their
intersection points $A_i$, $B_i$, $C_i$. Then the following holds:
\comment{\begin{align} 0= &\phantom{+}\;\;A_3A_4\cdot A_1C_1\cdot
A_2C_2\cdot B_1B_2\cdot B_3C_3\cdot B_4C_4
&-A_2A_4\cdot A_1C_1\cdot A_3C_3\cdot B_1B_3\cdot B_2C_2\cdot B_4C_4\notag\\
&+A_2A_3\cdot A_1C_1\cdot A_4C_4\cdot B_1B_4\cdot B_2C_2\cdot B_3C_3
&+A_1A_2\cdot A_3C_3\cdot A_4C_4\cdot B_1C_1\cdot B_2C_2\cdot B_3B_4\notag\\
&-A_1A_3\cdot A_2C_2\cdot A_4C_4\cdot B_1C_1\cdot B_2B_4\cdot B_3C_3
&+A_1A_4\cdot A_2C_2\cdot A_3C_3\cdot B_1C_1\cdot B_2B_3\cdot B_4C_4&.
\end{align}}
\begin{align}\label{eq:twelve}
&\;\;\phantom{+}
  A_1A_3\cdot A_2C_2\cdot A_4C_4\cdot B_1C_1\cdot B_2B_4\cdot B_3C_3\hskip-0.7em
&+\,A_2A_4\cdot A_1C_1\cdot A_3C_3\cdot B_1B_3\cdot B_2C_2\cdot B_4C_4\notag\\
&=A_3A_4\cdot A_1C_1\cdot A_2C_2\cdot B_1B_2\cdot B_3C_3\cdot B_4C_4\hskip-0.7em
&+\,A_2A_3\cdot A_1C_1\cdot A_4C_4\cdot B_1B_4\cdot B_2C_2\cdot B_3C_3\notag\\
&+A_1A_2\cdot A_3C_3\cdot A_4C_4\cdot B_1C_1\cdot B_2C_2\cdot B_3B_4\hskip-0.7em
&+\,A_1A_4\cdot A_2C_2\cdot A_3C_3\cdot B_1C_1\cdot B_2B_3\cdot B_4C_4&.
\end{align}
The segments are signed as before.
\end{theorem}
\begin{proof}
Add to the drawing the $x$ axis, which is neither parallel nor
orthogonal to any of the seven lines. For $i\in\{1,2,3,4\}$ let
$x_i,y_i,z_i$ be the projections of $A_i,B_i,C_i$ respectively onto
this axis. Then the theorem's claim~(\ref{eq:twelve}) can be
rewritten as follows:
\begin{scriptsize}\begin{align*}
&\;\;\phantom{+}
  (x_1 - x_3) (x_2 - z_2) (x_4 - z_4) (y_1 - z_1) (y_2 - y_4) (y_3 - z_3)\hskip-1em
&+(x_2 - x_4) (x_1 - z_1) (x_3 - z_3) (y_1 - y_3) (y_2 - z_2) (y_4 - z_4)\notag\\
&=(x_3 - x_4) (x_1 - z_1) (x_2 - z_2) (y_1 - y_2) (y_3 - z_3) (y_4 - z_4)\hskip-1em
&+(x_2 - x_3) (x_1 - z_1) (x_4 - z_4) (y_1 - y_4) (y_2 - z_2) (y_3 - z_3)\notag\\
&+(x_1 - x_2) (x_3 - z_3) (x_4 - z_4) (y_1 - z_1) (y_2 - z_2) (y_3 - y_4)\hskip-1em
&+(x_1 - x_4) (x_2 - z_2) (x_3 - z_3) (y_1 - z_1) (y_2 - y_3) (y_4 - z_4)&.
\end{align*}
\end{scriptsize}
With a help of a symbolic calculator it is easy to verify that the
above formula is equivalent~to:
\begin{equation}\label{eq:det}
\det\left(\begin{matrix}
x_1 - y_1 & x_1 - z_1 & z_1(x_1 - y_1) & y_1(x_1 - z_1) \\
x_2 - y_2 & x_2 - z_2 & z_2(x_2 - y_2) & y_2(x_2 - z_2)  \\
x_3 - y_3 & x_3 - z_3 & z_3(x_3 - y_3) & y_3(x_3 - z_3)  \\
x_4 - y_4 & x_4 - z_4 & z_4(x_4 - y_4) & y_4(x_4 - z_4)
\end{matrix}\right)=0.
\end{equation}
Observe that (\ref{eq:det}) is invariant under all
configuration-preserving permutations of the lines $l_{a,b,c}$ and
$l_{1,2,3,4}$, which is less than obvious given only the original
statement.

The proof consists of two substantially different cases.

\lead{Case I.} There exist two lines, say, $l_a$ and $l_b$, so that
$A_i\neq B_i$ for all $i\in \{1,2,3,4\}$. By appropriately scaling
and translating the $y$ axis we can make $l_a=\{x=y\}$. The rest of
the proof will be done in the homogenous coordinates. Let
$l_a=(1:-1:0)$ and $l_b=(\alpha:\beta:\gamma)$. Since $l_a$ and
$l_b$ are not equal, either $\gamma\neq 0$ or $\alpha\neq -\beta$.
For $i\in\{1,2,3,4\}$ the intersection point of $l_a$ and $l_i$
projects to $x_i$, and therefore $A_i=l_a\cap
(1:0:-x_i)=(x_i:x_i:1)$. Likewise, $B_i=(\beta y_i:-\gamma-\alpha
y_i:\beta)$. The two (distinct!) points uniquely define $l_i$:
$$
l_i=(x_i\beta+\gamma+\alpha y_i:-x_i\beta+\beta y_i:-\gamma
x_i-\alpha x_i y_i-\beta x_i y_i).
$$
The four lines $l_i$ for $i\in\{1,2,3,4\}$ intersect with the
vertical lines $(1:0:-z_i)$ at the following points that must be
collinear (as they all lie on $l_c$):
\begin{equation}\label{eq:coord}
C_i=(\beta z_i(x_i-y_i):\gamma(z_i-x_i)-\beta x_i(y_i-z_i)+\alpha
y_i(z_i-x_i):\beta(x_i-y_i)).
\end{equation}
The points are collinear if and only if the $3\times 4$ matrix whose
$i$th line is~(\ref{eq:coord}) has rank less than 3. Rank-preserving
transformation of these matrix reduce it to the matrix $M$ with the
following $i$th line (we may divide by $\beta\neq 0$ because $l_b$
is not vertical):
$$
z_i(x_i-y_i),(z_i-x_i)(\gamma+(\alpha+\beta) y_i), x_i-y_i.
$$
Since matrix $M$ has rank less than 3 for some $\alpha$, $\beta$, and
$\gamma$ subject to the condition that either $\gamma\neq 0$ or
$\alpha\neq -\beta$, it is equivalent to the following $4\times 4$
matrix being singular, which implies~(\ref{eq:det}):
$$
\left(\begin{matrix}
z_1(x_1-y_1) & z_1-x_1 & (z_1-x_1)y_1 & x_1-y_1\\
z_2(x_2-y_2) & z_2-x_2 & (z_2-x_2)y_2 & x_2-y_2\\
z_3(x_3-y_3) & z_3-x_3 & (z_3-x_3)y_3 & x_3-y_3\\
z_4(x_4-y_4) & z_4-x_4 & (z_4-x_4)y_4 & x_4-y_4\\
\end{matrix}\right).
$$

\lead{Case II.} For any two lines from $l_a,l_b,l_c$ there is some
$l_i$ going through their intersection. Reorder the lines so that
$A_1=C_1$, $B_2=C_2$, and $A_3=B_3$ (see Figure~\ref{fig:menelaus}).
\begin{figure}[htb]
\begin{center}
\psfrag{x11}{$A_1=C_1$}
\psfrag{x13}{$A_3=B_3$}
\psfrag{x14}{$A_4$}
\psfrag{x22}{$B_2=C_2$}
\psfrag{x24}{$B_4$}
\psfrag{x34}{$C_4$}
\psfrag{l1}{$l_1$}
\psfrag{l2}{$l_2$}
\psfrag{l3}{$l_3$}
\psfrag{l4}{$l_4$}
\psfrag{ll1}{$l_a$}
\psfrag{ll2}{$l_b$}
\psfrag{ll3}{$l_c$}
\includegraphics{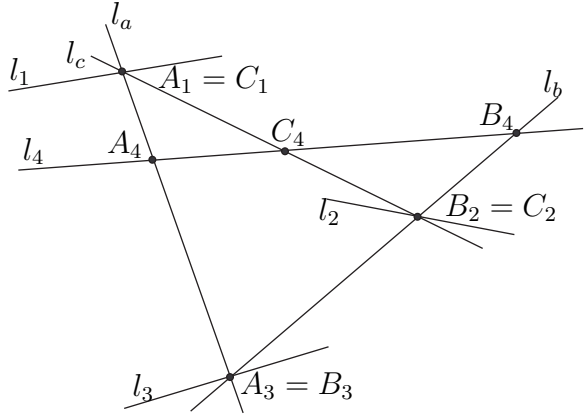}
\caption{Degenerate case.}\label{fig:menelaus}
\end{center}
\end{figure}

Plugging the identities into~(\ref{eq:twelve}) we reduce it to
$$
A_1A_3\cdot B_2B_4\cdot A_4C_4 = B_2B_3\cdot B_4C_4 \cdot A_1A_4,
$$
which is true by classic Menelaus' theorem.\QED\end{proof}

We note that Theorem~\ref{th:bipartite} is the ``minimal''
projective theorem that holds for the intersection points of two
sets of lines in general position. Indeed, it follows from the proof
that for any assignment of the eleven variables $x_{1,2,3,4}$,
$y_{1,2,3,4}$, $z_{1,2,3}$ there is a configuration of lines whose
intersection points project to those variables. This is because
there always exist $\alpha,\beta,\gamma$ such that $\gamma\neq 0$ or
$\alpha+\beta\neq 0$ and $C_1,C_2,C_3$ as defined
in~(\ref{eq:coord}) are collinear. \comment{
$$ \textrm{rank}\left(
\begin{matrix}
z_1(x_1-y_1) & (z_1-x_1)(\gamma+(\alpha+\beta)y_1) & x_1-y_1\\
z_2(x_2-y_2) & (z_2-x_2)(\gamma+(\alpha+\beta)y_2) & x_2-y_2\\
z_3(x_3-y_3) & (z_3-x_3)(\gamma+(\alpha+\beta)y_3) & x_3-y_3
\end{matrix}
\right)<3.
$$}
Hence no projection-invariant equation can be imposed on the lengths
of the segments that would not include all the twelve points. Other
configurations with as many or fewer intersection points are of no
use either: six lines intersecting two can project to any collection
of twelve points.

To highlight the projective nature of Theorem~\ref{th:bipartite}, we
may rewrite~(\ref{eq:twelve}) in a form that is invariant under
projection:
\begin{align*} 0=\phantom{+}
 &\frac{A_3A_4}{A_1A_4}\cdot \frac{B_1B_2}{B_1B_4}\cdot \frac{B_3C_3}{A_3C_3}\cdot \frac{B_4C_4}{A_4C_4} -
  \frac{A_2A_4}{A_1A_4}\cdot \frac{B_1B_3}{B_1B_4}\cdot \frac{B_2C_2}{A_2C_2}\cdot \frac{B_4C_4}{A_4C_4} +
  \frac{A_2A_3}{A_1A_4}\cdot \frac{B_2C_2}{A_2C_2}\cdot \frac{B_3C_3}{A_3C_3}\notag\\
+&\frac{A_1A_2}{A_1A_4}\cdot \frac{B_3B_4}{B_1B_4}\cdot
\frac{B_1C_1}{A_1C_1}\cdot \frac{B_2C_2}{A_2C_2}-
  \frac{A_1A_3}{A_1A_4}\cdot \frac{B_2B_4}{B_1B_4}\cdot \frac{B_1C_1}{A_1C_1}\cdot \frac{B_3C_3}{A_3C_3}+
  \frac{B_2B_3}{B_1B_4}\cdot \frac{B_1C_1}{A_1C_1}\cdot
  \frac{B_4C_4}{A_4C_4},
\end{align*}
where $\frac{AB}{CD}$ for parallel $AB$ and $CD$ equals
$\frac{|AB|}{|CD|}$ when the two segments have the same direction,
and $-\frac{|AB|}{|CD|}$ otherwise.

Finally, we transform~(\ref{eq:det}) to draw an analogy
with~(\ref{eq:menelausdet}):
\begin{equation}\label{eq:det-segments}
\det\left(\begin{matrix}
A_1B_1 & A_1C_1 & C_1C_1\cdot A_1B_1 & B_1B_1\cdot A_1C_1 \\
A_2B_2 & A_2C_2 & C_2C_1\cdot A_2B_2 & B_2B_1\cdot A_2C_2 \\
A_3B_3 & A_3C_3 & C_3C_1\cdot A_3B_3 & B_3B_1\cdot A_3C_3 \\
A_4B_4 & A_4C_4 & C_4C_1\cdot A_4B_4 & B_4B_1\cdot A_4C_4
\end{matrix}\right)=0.
\end{equation}

\end{document}